\documentclass[12pt]{article}
\usepackage{amsmath}

\usepackage{amsthm} 
\usepackage{amscd}
\theoremstyle{plain} 
\newtheorem{theorem}{Theorem}[section]
\newtheorem{corollary}[theorem]{Corollary}

\newtheorem{lemma}[theorem]{Lemma}

\numberwithin{equation}{section}

\newcommand{\al}{\alpha}

\newcommand{\sg}{\sigma}

\newcommand{\bq}{\mathbf {Q}}
\newcommand{\bz}{\mathbf {Z}}

\newcommand{\ot}{\otimes}

\begin{document}
\title{Leibniz Homology, Characteristic Classes and K-theory}
\author{Jerry M. Lodder\thanks{Supported
by the National Science Foundation, grant no. DMS-9704891.}}
\date{}
\maketitle

\section{Introduction}

In this paper we identify many striking elements in Leibniz
(co)homology which arise from characteristic classes and K-theory.
For a group $G$ and a field $k$ containing $\bq$, it is shown that all
primary characteristic classes, i.e. $H^*(BG; \, k)$, naturally inject
into certain Leibniz cohomology groups via an explicit chain map.
Moreover, if $f: A \to B$ is a homomorphism of algebras or rings, the
relative Leibniz homology groups $HL_*(f)$ are defined, and if in
addition $f$ is surjective with nilpotent kernel, $A$ and $B$ algebras
over $\bq$, then there is a natural surjection
$$  HL_{*+1}(gl(f)) \to HC_*(f),  $$
where $HC_*(f)$ denotes relative cyclic homology, and $gl(f): gl(A)
\to gl(B)$ is the induced map on matrices.
Here again, the above
surjection is realized via an explicit chain map, and offers a
relation between Leibniz homology and K-theory, since by work of
T. Goodwillie \cite{G}, there is an isomorphism 
$$  K_{*+1}(f) \ot \bq \to HC_*(f)  $$
between the relative theories when $f$ satisfies the above hypotheses.

Both explicit chain maps mentioned above involve an initial
homomorphism
$$  \varphi_* : HL_{*+1}(A) \to HH_*(A)  $$
from Leibniz to Hochschild homology, which on the chain level 
is simply a version of the
antisymmetrization map.  In many cases, $\varphi_*$ is seen to be
surjective, and from $HH_*(A)$ there are known maps to other types of
homologies, such as cyclic homology, or when $A=k[G]$, $HH_*(A)$ maps
to $H_*(BG)$.  The various homologies can be assembled into a curious
commutative diagram:
$$ \CD
HL_{*+1}(A) @>>>  HH_*(A)  \\
@VVV \ @VVV		  \\
H^{\rm{Lie}}_{*+1}(A) @>>> HC_*(A),
\endCD  $$
where $H^{\rm{Lie}}_*(A)$ denotes Lie-algebra homology.  

The map from Leibniz to Hochschild homology is studied in \S 2, while
the map to $H_*(BG)$ appears in \S 3, and the map to cyclic homology
in the relative case is in the final paragraph.  The appendix offers
an alternative calculation of the primitives of $HL_*(gl(A))$, and
fills a gap in a  proof in \cite{Lodder3}.  
Combined with the
author's previous results on foliations, Leibniz cohomology contains
primary characteristic classes (this paper), secondary characteristic
classes \cite{Lodder1}, and variations of secondary classes
\cite{Lodder2}.  In this way, the $HL_*$ theory can be viewed as a
``proto-homology.''

\section{From Leibniz to Hochschild Homology}

Recall that J.-L. Loday has defined Leibniz homology for the category
of Lie algebras, and more generally for the category of Leibniz
algebras \cite{LP}.  In this paper we begin with an associative
algebra $A$ over a commutative ring $k$, and consider $A$ as a Lie
algebra via its algebra structure:
$$  [a, \, b] = ab -ba, \ \ \ a, \ b \in A.  $$
In most cases $A$ will be unital, although that assumption is not
necessary for the definition of $HL_*(A)$, the Leibniz homology of
$A$.  The latter is defined as the homology of the chain complex
$CL_*(A)$:
$$ \CD
k @<0<< A @< [\ ,\ ]<< A^{\ot 2} @<<< \ldots @<<< A^{\ot n} @<d<< A^{\ot
(n+1)} @<<< \ldots \, ,
\endCD  $$
where
\begin{equation}
\begin{split}
& d(a_0 \ot a_1 \ot \, \ldots \, \ot a_n) = \\
& \sum_{0 \leq i < j \leq n} (-1)^{j+1}\, (a_0, \, a_1, \, \ldots, \,
a_{i-1}, \, [a_i, \, a_j], \, a_{i+1},\, \ldots , \, \hat{a}_j, \, \ldots, \, a_n).  
\end{split}
\end{equation}

The projection from the tensor powers of $A$ to the exterior powers
$$ A^{\ot n} \to A^{\wedge n} $$
induces a natural map
$$  HL_*(A) \to H^{\rm{Lie}}_*(A) $$
from Leibniz homology to Lie-algebra homology, where again the algebra
$A$ is viewed as a Lie algebra.  We now define a natural map
$$ HL_{*+1}(A) \to HH_*(A)  $$
to Hochschild homology, $HH_*(A)$, which when $k$ is a characteristic
zero field, yields a commutative diagram:
\begin{equation}  \label{1.1}
\CD
HL_{*+1}(A) @>>> HH_*(A)  \\
@VVV \ @VVV	          \\	
H^{\rm{Lie}}_{*+1}(A) @>>> HC_*(A),
\endCD 
\end{equation}
with $HC_*(A)$ denoting cyclic homology.

Recall that $HH_*(A)$ is the homology of the chain complex $CHH_*(A)$:
$$ \CD
A @<b<< A^{\ot 2} @<<< \, \ldots \, @<<< A^{\ot n} @<b<< A^{\ot (n+1)}
@<<< \, \ldots ,
\endCD  $$
where $b: A^{\ot (n+1)} \to A^n$ is defined using the face maps of the
cyclic bar construction \cite[1.1.1]{CH}
\begin{equation} \label{cyclicbar}
\begin{split}
& b = \sum_{i=0}^n (-1)^i \, d_i \\
& d_i(a_0, \, a_1, \, \ldots, \, a_n) = (a_0, \, a_1, \, \ldots, \,
      a_i a_{i+1}, \, \ldots, \, a_n), \ \ 0 \leq i < n \\
& d_n(a_0, \, a_1, \, \ldots, \, a_n) = (a_n a_0, \, a_1, \, a_2, \,
\ldots, \, a_{n-1}).
\end{split}
\end{equation}
Note that $b: A^{\ot 2} \to A$ is simply the bracket $[\ , \;
]:A^{\ot 2} \to A$, since
$$  b(a_0, \, a_1) = a_0 a_1 - a_1 a_0, $$
from which follows 
$$ HL_1(A) \simeq HH_0(A) \simeq HC_0(A) \simeq H^{\rm{Lie}}_1(A). $$
In fact, every arrow in diagram \eqref{1.1} is an isomorphism in the
special case $*=0$.

Consider now the map of chain complexes
\begin{equation}
\begin{split}
& \varphi : CL_{*+1}(A) \to CHH_*(A)  \\
& \varphi_n : A^{\ot (n+1)} \to A^{\ot (n+1)}  
\end{split}
\end{equation}
given by $\varphi_0 = \mathbf{1}$, $\varphi_1 =\mathbf{1}$, and for $n
\geq 2$,
\begin{equation} \label{1.2}
\begin{split}
  & \varphi_n(a_0, \, a_1, \, \ldots, \, a_n) = \\
  & \sum_{\sigma \in S_n} ({\rm{sgn}} \, \sigma) \, (a_0, \,
  a_{\sigma^{-1}(1)}, \, a_{\sigma^{-1}(2)}, \, \ldots, \,
  a_{\sg^{-1}(n)}), 
\end{split}
\end{equation}
where $S_n$ is the symmetric group on $n$ letters.  Of course, formula
\eqref{1.2} could be used to define $\varphi_1$, since $S_1$ is the
trivial group.  Also, the summation remains invariant if $\sigma^{-1}$
is replaced with $\sigma$ in all subscripts, which reconciles various
descriptions of antisymmetrization maps in the literature.
\begin{lemma} \label{1.3}
The $k$-module homomorphism
$$  \varphi : CL_{*+1}(A) \to CHH_*(A)  $$
is a map of chain complexes.  
\end{lemma}
\begin{proof}
It follows at once that $b \varphi_1 = \varphi_0 d$.  For $n \geq 2$,
\begin{equation}
\begin{split}
  & b \circ \varphi_n (a_0, \, a_1, \, a_2, \, \ldots, \, a_n) = \\
  & \ \ \ \ \ \sum_{\sigma \in S_n} ({\rm{sgn}} \, \sigma) \, \Big( (a_0
  a_{\sigma^{-1}(1)}, \, a_{\sigma^{-1}(2)}, \, \ldots \,
  a_{\sg^{-1} (n)})\\
  &\ \ \ \ \  -(a_0, \, a_{\sg^{-1}(1)}a_{\sg^{-1}(2)}, \, a_{\sg^{-1}(3)}, \,
  \ldots, \, a_{\sg^{-1}(n)})  \\
  &\ \ \ \ \   + \cdots + (-1)^i (a_0, \, \ldots, \, a_{\sg^{-1}(i)}
  a_{\sg^{-1}(i+1)}, \, \ldots, \, a_{\sg^{-1}(n)})  \\
  &\ \ \ \ \   + \cdots + (-1)^n ( a_{\sg^{-1}(n)}a_0, \, a_{\sg^{-1}(1)}, \,
  \ldots, \, a_{\sg^{-1}(n-1)}) \Big).
\end{split}
\end{equation}
On the other hand,
\begin{equation}
\begin{split}
  & \varphi_{n-1}\circ d (a_0, \, a_1, \, \ldots, \, a_n) =  \\
  & \ \ \ \ \ \varphi_{n-1}\Big( \sum_{j=1}^n (-1)^{j+1} \, ( [a_0, \,
  a_j], \, a_1, \, \ldots, \, \hat{a}_j, \, \ldots, \, a_n)\Big)  \\
  & \ \ \ \ \ + \varphi_{n-1}\Big( \sum_{1 \leq i < j \leq n}
  (-1)^{j+1} \, (a_0, \, a_1, \, \ldots, \, [a_i, \, a_j], \ldots, \,
  \hat{a}_j, \, \ldots, \, a_n) \Big).
\end{split}
\end{equation}
Note that the terms
\begin{equation}
\begin{split}
\sum_{\sigma \in S_n} & ({\rm{sgn}} \, \sigma) \, \Big( (a_0
  a_{\sigma^{-1}(1)}, \, a_{\sigma^{-1}(2)}, \, \ldots \, a_{\sg^{-1}(n)})\\
  & + (-1)^n ( a_{\sg^{-1}(n)}a_0, \, a_{\sg^{-1}(1)}, \,
  \ldots, \, a_{\sg^{-1}(n-1)}) \Big)
\end{split}
\end{equation}
match the terms
$$ \varphi_{n-1}\Big( \sum_{j=1}^n (-1)^{j+1} \, ( [a_0, \,
  a_j], \, a_1, \, \ldots, \, \hat{a}_j, \, \ldots, \, a_n)\Big),$$
and the remaining terms in the sums for $b \circ \varphi_n$ and
$\varphi_{n-1} \circ d$ also match.  Thus, 
$$  b \circ \varphi_n = \varphi_{n-1} \circ d .  $$
\end{proof}
\begin{corollary} \label{1.4}
There is a natural induced homomorphism
$$  \varphi_* : HL_{*+1}(A) \to HH_*(A) .  $$
\end{corollary}

The chain map $\varphi$ belongs to the genre of constructions known as
antisymmetrization maps \cite[1.3.4]{CH}.  In fact, it follows from
\cite[1.3.5]{CH} that $\varphi$ descends to the exterior powers on the
domain
$$ \epsilon_n : A \ot \Lambda^n (A) \to A^{\ot (n+1)},  $$
and there is an induced map
\begin{equation} \label{1.5}
  \epsilon_* : H^{\rm{Lie}}_*(A; \, A) \to HH_*(A)
\end{equation}
where $H^{\rm{Lie}}_*(A; \, A)$ denotes Lie-algebra homology with
coefficients in the adjoint representation:
$$  {\rm{ad}}(a)(b) = [a, \, b], \ \ a, \ b \in A.  $$
For completeness, recall that $H^{\rm{Lie}}_*(A; \, A)$ is the homology of
the complex
\begin{align*}
  & A \overset{[\; , \;]}{\longleftarrow} A \ot A  \longleftarrow A \ot A^{\wedge 2}
  \longleftarrow \, \ldots \,  A \ot A^{\wedge (n-1)}
  \overset{d}{\longleftarrow} A \ot A^{\wedge n} \longleftarrow \,
  \ldots , \\
  & d (a_0 \ot a_1 \wedge a_2 \wedge \, \ldots \, a_n) = \\ 
  &  \ \ \ \ \ \sum_{1 \leq i
  < j \leq n} (-1)^{j+1} \, (a_0 \ot a_1 \wedge \, \ldots \, a_{i-1}
   \wedge [a_i, \, a_j] \wedge a_{i+1} \wedge \, \ldots \, \hat{a}_j \,
   \ldots \, \wedge a_n)  \\
  & \ \ \ \ \ + \sum_{j=1}^n (-1)^{j+1} \, \big( {\rm{ad}}(a_0)(a_j) \ot
  a_1 \wedge a_2 \wedge \, \ldots \, \hat{a}_j \, \ldots \wedge
  a_n\big) .
\end{align*}

The projection $A^{\ot (n+1)} \to A \ot A^{\wedge n}$ induces a map on
homology
\begin{equation} \label{1.6}
  HL_{*+1} (A) \to H^{\rm{Lie}}_*(A; \, A).
\end{equation}
\begin{lemma} \label{1.7}
There is a commutative diagram
$$ \CD
  HL_{*+1}(A)  \\
  @VVV \hskip-1in \searrow\varphi_*   \\
  H^{\rm{Lie}}_* (A; \, A) @>>\epsilon_*> HH_*(A) 
\endCD  $$  
\end{lemma}
\begin{proof}
This follow from corollary \eqref{1.4} and maps \eqref{1.5} and
\eqref{1.6}
\end{proof}

When $k$ is a characteristic zero field, the cyclic homology of $A$
may be computed from the complex $C_*^{\lambda}(A)$ \cite[2.1.4]{CH}:
$$  A \overset{b}{\longleftarrow} A^{\ot 2}/(1-t) \longleftarrow \,
    \ldots \, A^{\ot n}/(1-t) \overset{b}{\longleftarrow} A^{\ot
    (n+1)}/(1-t) \longleftarrow \, \ldots ,$$
where $\bz /(n+1)$ acts on $A^{\ot (n+1)}$ via
$$  t(a_0, \, a_1, \, \ldots, \, a_n) = (-1)^n (a_n, \, a_0, \, a_1,
    \, \ldots, \, a_{n-1}).  $$
There is a chain map \cite[10.2.3]{CH}
\begin{equation} \label{1.8}
\begin{split}
  & \theta : \Lambda^{*+1} (A) \to C^{\lambda}_* (A)  \\
  & \theta : \Lambda^{n+1} (A) \to A^{\ot (n+1)}/(1-t)  \\
  & \theta(a_0 \wedge a_1 \wedge \, \ldots \, \wedge a_n) = \sum_{\sigma \in
  S_n} ({\rm{sgn}} \, \sigma) (a_0, \, a_{\sigma^{-1}(1)}, \,
  a_{\sg^{-1}(2)}\, \ldots , \, a_{\sigma^{-1}(n)}),
\end{split}
\end{equation}
which induces a homomorphism
$$  \theta_* : H^{\rm{Lie}}_{*+1}(A) \to HC_*(A).  $$
Recall that there is a natural map $I : HH_*(A) \to HC_*(A)$
\cite[2.2.1]{CH}, which in the characteristic zero case is induced by
the projection
\begin{equation}
\begin{split}
  & CHH_*(A) \to C_*^{\lambda}(A)  \\
  & A^{\otimes (n+1)} \to A^{\otimes (n+1)}/(1-t).
\end{split}
\end{equation}
\begin{lemma} \label{1.9}
When $k$ is a characteristic zero field, there is a commutative
diagram
$$ \CD
   HL_{*+1}(A) @>\varphi_*>>  HH_*(A)  \\
   @V\pi VV  \  @VIVV  \\
   H^{\rm{Lie}}_{*+1}(A) @>\theta_*>>  HC_*(A).
\endCD  $$
\end{lemma}
\begin{proof}
This follows immediately, since both $\varphi_*$ and $\theta_*$ are
induced by antisymmetrization maps, and both $\pi$ and $I$ are induced
by projections.  
\end{proof}
Of course, when $*=0$, every arrow in lemma \eqref{1.9} is an
isomorphism.

The map $\varphi_* : HL_{*+1}(A) \to HH_*(A)$ has a nice
interpretation when $A$ is a smooth algebra over a Noetherian ring
$k$.  By definition of smooth, $A$ is assumed to be commutative and
unital over $k$ (See \cite[3.4.1]{CH} for a discussion of smooth algebras.)
Let $\Omega^*_{A|k}$ be the $A$-module of differential forms, where
$$  \Omega^0_{A|k} = A, \ \ \ \Omega^n_{A|k} = \Lambda^n_A
(\Omega^1_{A|k}),  $$
and $\Omega^1_{A|k}$ form the K\"ahler differentials, generated as an
$A$-module by the symbols $da$, $a \in A$, subject to the relations
\begin{equation}
\begin{split}
  & d(\lambda a + \mu b) = \lambda \, da + \mu \, db, \ \ \lambda, \,
  \mu \in k, \ \ a, \, b \in A  \\
  &  d(ab) = a(db) + b (da), \ \ a, \, b \in A 
\end{split}
\end{equation}
By a theorem of Hochschild, Kostant, and Rosenberg \cite{HKR}, when
$A$ is smooth over $k$, the antisymmetrization 
\begin{equation} \label{1.10}
  \epsilon_* : \Omega^*_{A|k} \to HH_*(A) 
\end{equation}
is an isomorphism of graded algebras.  Of course, for any commutative
algebra $A$, $[a, \, b] = 0$ for all $a, \; b \in A$, and
\begin{equation}  \label{1.11}
  HL_*(A) \simeq T(A),
\end{equation}
where $T(A) = \sum_{k  \geq 0} A^{\ot k}$ denotes the tensor algebra
on $A$.  
\begin{lemma}  \label{1.12}
If $A$ is smooth over $k$, then the natural map
$$  \varphi_*:HL_{*+1}(A) \to HH_*(A)  $$
is surjective.
\end{lemma}
\begin{proof}
There is a surjective $k$-module homomorphism
$$  HL_{n+1}(A) \simeq A^{\ot (n+1)} \to \Omega^n_{A|k}  $$
given on homogeneous elements by
$$  p(a_0 \ot a_1 \ot \, \ldots \, \ot a_n) = a_0 \, da_1 \wedge da_2
\wedge \, \ldots \, \wedge da_n .  $$
The result now follows from the commutativity of the diagram
$$ \CD
  HL_{*+1}(A)  \\
  @VpVV \hskip-1in \searrow \, \varphi_*  \\
  \Omega^*_{A|k} @>>\epsilon_*>  HH_*(A).
\endCD  $$
\end{proof}
When $A$ is the smooth algebra $C^{\infty}(M)$ of real-valued
differentiable functions on a differentiable manifold $M$, the reader
is invited to identify the various homology groups in lemma
\eqref{1.9} using the calculations of $HC_*( C^{\infty}(M))$ in terms
of de Rham cohomology \cite[3.4.12]{CH}.  

In general, the map $\varphi_* : HL_{*+1}(A) \to HH_*(A)$ is not
surjective, even in the commutative case, since from the definition of $\varphi$
$$  {\rm{Im}} \, [\varphi_*(HL_{n+1}(A))] \subset HH^{(n)}_n(A), $$
where $HH^{(n)}_n(A)$ denotes the n-th summand in the
$\lambda$-decomposition of Hoch-schild homology over $\bq$.  See
\cite[4.5]{CH} for details about this decomposition.  Of course,
$HH^{(n)}_n(A)$ is not necessarily isomorphic to $HH_n(A)$.  In the
next section we investigate another important case in which
$\varphi_*$ is surjective.

\section{Relation to Characteristic Classes}

In this paragraph we prove that the natural map
$$ \varphi_* : HL_{*+1}(gl(A)) \to HH_*(gl(A)) \simeq HH_*(A)  $$
is onto for a unital algebra over a characteristic zero field $k$.
Furthermore, when $A$ is the group ring $k[G]$, there is a surjective
homomorphism
$$  HH_*( k[G]) \to H_*(k[G]) := H_*(BG; \, k),  $$
and, in fact, $H_*(BG; \, k)$ is a direct summand of $HH_*(k[G])$.
Thus, on cohomology
$$  H^*(BG; \, k) \to HL^{*+1}(gl(A))  $$
is injective, and $HL^{*+1}(gl(A))$ contains all characteristic
classes in $$ H^*(BG; \, k). $$  

Specifically, let
$$ gl(A) = \underset{\overset{\longrightarrow}{n}}{\rm{lim}} \, gl_n(A) $$ 
be the Lie algebra of infinite matrices over $A$ with finitely many
nonzero entries.  Note that the Lie algebra structure on $gl(A)$ is
actually induced from the ring structure of
$$  M(A) = \underset{\overset{\longrightarrow}{n}}{\rm{lim}} \, M_n(A), $$
where $M_n(A) = gl_n(A)$ is the ring of all $n \times n$ matrices over
$A$.  Recall that \cite{Cu} \cite[10.6.5]{CH} 
\begin{equation} \label{2.1}
\begin{split}
  HL_*(gl(A)) \simeq T(HH_*(A)[1]),  
\end{split}
\end{equation}
where $T$ denotes the tensor algebra.  See also \cite{Lodder3}.
To understand the map $\varphi_*$ on homology, the isomorphism in
equation \eqref{2.1} must be understood on the chain level.  The chain
complex $CL_*(gl(A))$ is quasi-isomorphic to $L_*(A)$ \cite[10.6.7]{CH}:
\begin{equation} \label{2.2}
\begin{split}
  k[S_1] \underset{k}{\ot} A \longleftarrow k[S_2]\underset{k}{\ot} A^{\ot
  2} \longleftarrow \, \ldots \, \longleftarrow k[S_{n}] \underset{k}{\ot}
  A^{\ot n} \longleftarrow \, \ldots ,
\end{split}
\end{equation}
where $S_n$ denotes the symmetric group on $n$ letters.  There is an
explicit chain map
\begin{equation}
\begin{split}
  & CL_*(gl(A)) \to L_*(A)  \\
  & [gl(A)]^{\ot n} \to \big( [gl(A)]^{\ot n} \big)_{gl(k)}
  \overset{\Theta}{\rightarrow} k[S_n]\ot A^{\ot n}
\end{split}
\end{equation}
with $[gl(A)]^{\ot n} \to \big( [gl(A)]^{\ot n} \big)_{gl(k)}$ being the
projection onto the quotient by the adjoint action.  For $\sigma \in
S_n$ and $E_{ij}^a$ the elementary matrix with only one possible nonzero entry
$a \in A$ in the $ij$ position,
\begin{equation} \label{2.3}
\begin{split}
  \Theta(E^{a_1}_{1 \, \sg (1)} \ot E^{a_2}_{2 \, \sg (2)} \ot \,
  \ldots \, \ot E^{a_n}_{n \, \sg (n)}) = \sg \ot (a_1, \, a_2, \, \ldots,
  \, a_n)
\end{split}
\end{equation}
For more details, see \cite[10.2.11]{CH}.

Now, $HH_*(A)[1]$ is the direct summand of $HL_*(gl(A))$ that arises
from the homology of the complex $P_* (A)$:
\begin{equation}  \label{2.4}
\begin{split}
k[U_1] \underset{k}{\ot} A \longleftarrow k[U_2]\underset{k}{\ot} A^{\ot
  2} \longleftarrow \, \ldots \, \longleftarrow k[U_{n}] \underset{k}{\ot}
  A^{\ot n} \longleftarrow \, \ldots ,
\end{split}
\end{equation}
where $U_n$ is the conjugacy class of the cyclic shift in $S_n$.  
The complex $P_*(A)$ is a summand of $L_*(A)$
\cite{Lodder3}, and an alternative calculation of $H_*(P_*(A))$
appears in the appendix.  Also needed for the identification of Im$(\varphi_*)$
is an explicit description of the trace map isomorphism from Morita invariance
$$  {\rm{tr}}_* : HH_*(gl(A)) \to HH_*(A) . $$
For ${}_pM = [{}_pm_{ij}] \in gl(A),$
\begin{equation}
\begin{split}
  & {\rm{tr}} ( {}_0M \ot {}_1M \ot \, \ldots \, \ot {}_nM) = \\
  & \ \ \ \ \ \sum {}_0m_{i_1 \, i_2} \ot {}_1m_{i_2 \, i_3} \ot {}_2
  m_{i_3 \, i_4} \ot
  \, \ldots \, \ot {}_n m_{i_{n+1}\, i_1} \, , 
\end{split}
\end{equation}
where the sum is over all indices $i_1$, $i_2$, $\ldots,$ $i_{n+1}$.
\begin{lemma}
Let $A$ be a unital algebra over a characteristic zero field $k$.
Then the composition
$$  {\rm{tr}}_* \circ \varphi_* : HL_{*+1}(gl(A)) \to HH_*(A)  $$
is surjective.
\end{lemma}
\begin{proof}
The summand of $HL_{*+1}(gl(A))$ isomorphic to $HH_*(A)[1]$ can be
represented via chains which are $k$-linear combinations of terms:
$$  E^{a_1}_{1 \, \tau(1)} \ot  E^{a_2}_{2 \, \tau(2)} \ot \, \ldots
	\, \ot  E^{a_n}_{n \, \tau(n)} ,$$
where $\tau$ is the cyclic shift given by the cycle 
$(1, \, 2, \, 3, \, \ldots, \, n)$ \cite{Lodder3}.  It is enough to
compute
\begin{align*}
  & {\rm{tr}} \circ \varphi (E^{a_1}_{1 \, \tau(1)} \ot  E^{a_2}_{2
	\, \tau(2)} \ot \, \ldots \, \ot  E^{a_n}_{n \, \tau(n)}) \\
  & = {\rm{tr}} \circ \varphi (E^{a_1}_{1 \, 2} \ot  E^{a_2}_{2 \, 3} \ot \, \ldots
	\, \ot  E^{a_n}_{n \, 1})  \\
  & = a_1 \ot a_2 \ot \, \ldots \, \ot a_n .
\end{align*}
Thus, any element in $HH_*(A)$ can be represented as a chain in
$CL_{*+1}(gl(A))$.
\end{proof}

Consider the case of a group ring $A = k[G]$, and the associated bar
construction $\{ G^n \}_{n \geq 0}$.  The face maps $d_i : G^n \to
G^{n-1}$ are given by
$$  d_i(g_1, \, g_2, \, \ldots, \, g_n) =\begin{cases}
	(g_2, \, g_3, \, \ldots, \, g_n) & i = 0, \\
	(g_1, \, \ldots, \, g_i g_{i+1}, \, \ldots, \, g_n) & 1 \leq i
			\leq n-1,  \\
	(g_1, \, g_2, \, \ldots, \, g_{n-1}) & i =n ,
	\end{cases}  $$
with $H_*(BG)$ denoting the homology of the complex
$$  B_*(G) = k[G^n] \simeq (k[G])^{\ot n}, \ \ \ n \geq 0.  $$
There are natural simplicial maps
\begin{align*}
  & \pi : CHH_*(k[G]) \to B_*(G)  \\
  & \pi : k[G^{n+1}] \to k[G^n]  \\
  & \pi (g_0, \, g_1, \, g_2, \, \ldots, \, g_n) = (g_1, \, g_2, \,
  \ldots, \, g_n) 
  &  \\
  & \iota : B_*(G) \to CHH_*(k[G])  \\
  & \iota : k[G^n] \to k[G^{n+1}]  \\
  & \iota (g_1, \, g_2, \, \ldots, \, g_n) = ( (g_1 g_2 \, \cdots \,
    g_n)^{-1}, \, g_1, \, g_2, \, \ldots, \, g_n),
\end{align*}
and the composition $\pi \circ \iota$ is the identity on $B_*(G)$.
Thus follows the known lemma \cite{B}
\begin{lemma} \cite{B}  The graded group $H_*(BG)$ is a direct summand of
$$HH_*(k[G])$$ for any coefficient ring $k$.
\end{lemma}

By composing the chain maps $\pi$, tr, and $\varphi$, we have
\begin{theorem} \label{2.7} Let $k$ be a characteristic zero field and
$A = k[G]$.  Then
$$  (\pi \circ {\rm{tr}} \circ \varphi)_* : HL_{*+1}(gl(A)) \to
	H_*(BG; \, k)  $$
is surjective, and
$$ (\pi \circ {\rm{tr}} \circ \varphi)^* : H^*(BG; \, k) \to
	HL^{*+1}(gl(A))  $$
in injective.  
\end{theorem}
Thus, all characteristic classes over $k$ naturally inject into
Leibniz cohomology.

\section{Relation to K-theory}

A fundamental theorem of T. Goodwillie \cite{G} relates relative
algebraic K-theory, $K_*(f)$, to relative cyclic homology, $HC_*(f)$.
In particular, if $f: R \to S$ is a homomorphism of simplicial rings
such that the induced map $\pi_0(R) \to \pi_0(S)$ is surjective with
nilpotent kernel, then 
\begin{equation} \label{3.1}
  K_n(f) \ot \bq \simeq HC_{n-1}(f) \ot \bq .
\end{equation}
Although for Leibniz homology we are working in the category of
discrete rings, the following results may be extended to the
simplicial case.  For a homomorphism $f : A \to B$ of discrete rings,
we define the relative Leibniz homology groups, $HL_*(f)$, we then consider
the map on matrices $gl(f) : gl(A) \to gl(B)$, and show that over a
characteristic zero field, the composition
$$  ({\rm{tr}}\circ \varphi)_* : HL_{*+1}(gl(f)) \to HH_*(f)  $$
is onto.  If furthermore, $f: A \to B$ is surjective with nilpotent
kernel, then
$$  I : HH_*(f) \to HC_*(f)  $$
is also surjective, as well as the composition
\begin{equation} \label{3.2}								
  I \circ ({\rm{tr}} \circ \varphi)_* : HL_{*+1}( gl(f)) \to HC_*(f).
\end{equation}
Above, all relative homology groups are taken with coefficients in a
characteristic zero field.

Recall briefly the construction of relative homology in a
general setting \cite[p. 46--47]{M}.  Let $f : C_n \to C'_n$, $n \geq
0$, be a homomorphism between any two chain complexes $C_*$ and
$C'_*$, and define the mapping cone of $f$ as
$$  M_n(f) = C_{n-1} \oplus C'_n, \ \ \ n \geq 1 ,  $$
with boundary map 
$$  \partial (c, c') = ( - \partial c, \, \partial c' + fc )  $$
By definition, $H_*(f)$ are the homology groups of $M_*[-1] =M_{*+1}$,
which fit into a long exact sequence
\begin{equation} \label{3.3}
\begin{split}
  & \ldots \, \longrightarrow H_n(C_*) \overset{f_*}{\longrightarrow} H_n(C'_*)
  \overset {\al_*}{\longrightarrow} H_n(f)
  \overset{p_*}{\longrightarrow} H_{n-1}(C_*) \longrightarrow \,
  \ldots \\
  & \ldots \, \longrightarrow H_1 (f)
  \overset{p_*}{\longrightarrow} H_0(C) \overset{f_*}{\longrightarrow}
  H_0(C'_*) \overset{\al_*}{\longrightarrow} H_0(f) \longrightarrow 0,
\end{split}
\end{equation}
where $p: M_*(f) \to C_{*-1}$ and $\al : C'_* \to M_*(f)$ are given by
$$  p((c, \, c')) = c, \ \ \ \al(c') = (0, \, c').  $$

If $f: A \to B$ is a homomorphism of discrete rings, then the above
construction yields $HL_*(f)$ by considering the chain map $f: CL_*(A)
\to CL_*(B)$.  Moreover, the relative homology group construction is
functorial, and there are natural maps
\begin{align*}
  &  \varphi_* : HL_{*+1}(f) \to HH_*(f)  \\
  & {\rm{tr}}_* : HH_*(gl(f)) \to HH_*(f) \\
  & I : HH_*(f) \to HC_*(f).
\end{align*}
\begin{theorem} Suppose that $A$ and $B$ are unital algebras over a
  characteristic zero field $k$, and $f: A \to B$ is an algebra
  homomorphism.  For the map $gl(f) : gl(A) \to gl(B)$, the
  composition
$$  ({\rm{tr}} \circ \varphi)_* : HL_{*+1}( gl(f)) \to HH_*( gl(f))
	\to HH_*(f)  $$
is surjective.
\end{theorem}
\begin{proof}
From equation \eqref{2.1} 
\begin{align*}
  & HL_*( gl(A)) \simeq T( HH_*(A) [1]) \\
  & HL_*( gl(B)) \simeq T( HH_*(B) [1])
\end{align*}
From equation \eqref{2.4} the summand of $HL_*(gl(A))$ isomorphic to
$HH_*(A)[1]$ can be computed from the complex $P_*(A)$, and similarly
for $gl(B)$.  Consider the $k$-linear homomorphism
\begin{align*}
  & P(f) : P_*(A) \to P_*(B)  \\
  & P(f) : k[U_n] \ot A^{\ot n} \to k[U_n] \ot B^{\ot n} \\
  & P(f) ( \sg \ot (a_1, \, a_2, \, \ldots, \, a_n)) = \sg \ot (
  f(a_1), \, f(a_2), \, \ldots, \, f(a_n)).  
\end{align*}
We have the following commutative diagram with exact rows:
$$ \CD
  \longrightarrow H_{n+1}(P_*(A)) @>>> H_{n+1}(P_*(B)) @> \al_* >> 
	H_{n+1}(P(f)) \longrightarrow  \\
   @VVV \ @VVV \ @VVV \ \\
  \longrightarrow HH_n(A) @>>> HH_n(B) @> \al_* >> HH_n(f)
  \longrightarrow  \, ,
\endCD  $$ 
where the dimensions in the top row are inherited from the chain
complex for Leibniz homology.  
By the 5-lemma, $H_*(P(f)) \simeq HH_*(f)[1]$.  Let $z$ denote an
element in the mapping cone 
$$  CL_*(gl(A)) \to CL_*(gl(B))  $$
of the form
$$ ( E^{a_1}_{1\, \tau(1)} \ot E^{a_2}_{2\, \tau(2)} \ot \, \ldots \,
   \ot E^{a_n}_{n\, \tau(n)}, \ E^{b_1}_{1\, \tau'(1)} \ot E^{b_2}_{2
   \, \tau'(2)} \ot \, \ldots \, \ot E^{b_{n+1}}_{n+1, \, \tau'(n+1)}),$$
where $\tau$ is the cyclic shift in $S_n$ and $\tau'$ is the cyclic
shift in $S_{n+1}$.  The chain map $\Theta$ in equation \eqref{2.3}
may be defined on the respective mapping cones, and
$$ \Theta (z) = \big( \tau \ot (a_1, \, a_2, \, \ldots, \, a_n), \ \tau'
	\ot (b_1, \, b_2, \, \ldots, \, b_{n+1}) \big),  $$
where the latter is in fact an element of $M(P(f))$.  Also, at the
level of mapping cones, we have
$$ ({\rm{tr}} \circ \varphi)(z) = \big( (a_1, \, a_2, \, \ldots, \, a_n), \
	(b_1, \, b_2, \, \ldots, \, b_{n+1}) \big) .  $$
The theorem follows, since the isomorphism
$$  H_*(P(f)) \simeq HH_*(f)[1]  $$
is realized by sending the class of $\Theta (z)$ to the class of
$({\rm{tr}}\circ \varphi)(z)$.
\end{proof}

For relative cyclic homology, there is a long exact sequence
$$ \ldots \, \longrightarrow HH_n(f) \overset{I}{\longrightarrow}
   HC_n(f) \overset{S}{\longrightarrow} HC_{n-2}(f)
   \overset{B}{\longrightarrow} HH_{n-1}(f) \overset{I}{\longrightarrow}
   \, \ldots  \, .  $$
By work of Goodwillie, we have
\begin{lemma} \label{3.5}
If $f : A \to B$ is a homomorphism of unital algebras over a
characteristic zero field with nilpotent kernel, then the map
$$  I : HH_n (f) \to HC_n (f)  $$ 
is surjective.
\end{lemma}
\begin{proof}
From \cite[p. 399]{G}, the map $S: HC_n (f) \to HC_{n-2} (f)$ is
zero.  
\end{proof}
\begin{corollary} \label{3.6}
Under the hypotheses of lemma \eqref{3.5}, the natural map
$$  I \circ ({\rm{tr}} \circ \varphi)_* : HL_{*+1}(gl(f)) \to HC_*(f)  $$
is surjective.
\end{corollary}

\section{Appendix:  The Homology of $P_*(A)$}

The results of this paper concerning the surjectivity of the
homomorphism
$$  HL_{*+1}(gl(A)) \to HH_*(A)  $$
rely heavily on the calculation of the homology of the complex $P_*(A)$:
$$ k[U_1] \underset{k}{\ot} A \longleftarrow k[U_2]\underset{k}{\ot} A^{\ot
  2} \longleftarrow \, \ldots \, \longleftarrow k[U_{n+1}] \underset{k}{\ot}
  A^{\ot (n+1)} \longleftarrow \, \ldots , $$
which form the primitive elements of $HL_*(gl(A))$.  In this appendix
we offer a calculation of $H_*(P_*(A))$ which fills a gap in a
previous proof, namely lemma (2.6) of \cite{Lodder3}.  We prove that
if $A$ is a unital $k$-algebra, then
$$  H_*(P_*(A)) \simeq HH_*(A),  $$
where now $P_n(A) = k[U_{n+1}] \underset{k}{\ot} A^{\ot (n+1)}$.

Recall that $U_n$ is the conjugacy class of the cyclic shift in the
symmetric group $S_n$.  The collection $\{ U_{n+1} \}_{n \geq 0}$ form
a presimplicial set with face maps $d_i$ \cite{Lodder3}, but lack
degeneracies.  Let $N^{\rm cy}(A)$ be the cyclic bar construction on
$A$ with faces given in equation \eqref{cyclicbar}.  Let
$$ k[U_{*+1}] \ot N^{\rm cy}(A)  $$
be the presimplicial $k$-module with
$$ \{ \, k[U_{*+1}] \ot N^{\rm cy} (A) \, \}_n = k[U_{n+1}]\ot A^{\ot
(n+1)},  $$
and face maps 
$$  d_i ( \sg \ot \vec{a}) = d_i(\sg ) \ot d_i (\vec{a}).  $$
The complex $P_*(A)$ is simply $k[U_{*+1}] \ot N^{\rm cy}(A)$ together
with its boundary map constructed as the alternating sum of the face
maps.  In \cite{Lodder3} it is shown that the complex $k[U_{*+1}]$ is
acyclic, but without degeneracies, the Eilenberg-Zilber and K\"unneth theorems cannot
be applied to calculate $H_*(P_*(A))$.  To remedy this, we invoke the
Dold-Kan functor \cite{Lamotke}, which to any chain complex $K_*$
(over $k$), associates the simplicial $k$-module
$$  D_n(K_* ) = {\rm Hom} ( C(\Delta [n]), \ K_*),  $$
where Hom denotes chain maps over $k$.  Also, $\Delta [n]$ is the
simplicial model for the $n$-simplex, $C'( \Delta [n])$ is the free
$k$-module on the elements of $\Delta [n]$, and
$$  C(\Delta [n]) = C'( \Delta [n])/\, ^eC(\Delta [n]),  $$
with $^eC(\Delta [n])$ denoting the submodule of $C'( \Delta [n])$
generated by the degenerate elements.  See \cite{Lamotke} for further
details.  

Using properties of the Dold-Kan functor, the complexes
$$  k[U_{*+1}] \ot N^{\rm cy}(A), \ \ \ D_*(k[U_{*+1}]) \ot N^{\rm cy}(A)$$
are quasi-isomorphic, as well as the complexes
$$  k[U_{*+1}] \ot N^{\rm cy}(A), \ \ \ D_*\big( k[U_{*+1}]\ot N^{\rm
cy}(A) \big).  $$
Define a chain map
\begin{equation*}
\begin{split}
& \psi : D_*(k[U_{*+1}]) \ot N^{\rm cy}(A) \to D_*\big( k[U_{*+1}]\ot N^{\rm
cy}(A) \big) \\
& \psi : {\rm{Hom}} (C(\Delta [n]), \; k[U_{*+1}]) \ot A^{\ot (n+1)} \to
{\rm{Hom}}( C(\Delta [n]), \; k[U_{*+1}]\ot A^{*+1})
\end{split}
\end{equation*}
as follows.  Recall that $\delta^i : [n-1] \to [n]$, $i = 0$, 1, 2,
$\ldots$, $n$ is given by
$$  \delta^i (j) = \begin{cases} j, & j < i \\
 				 j+1, & j \geq i. \end{cases}$$
Let $\al \in (\Delta [n])_q$, $\, \al: [q] \to [n]$, $\, q \leq n$.  Then
\begin{equation*}
\begin{split}
& \al = \delta^{i_{n-1}} \circ \, \ldots \, \circ \delta^{i_2} \circ
\delta^{i_1} \\
& [q] \overset{\delta^{i_1}}{\longrightarrow} [q+1]
\overset{\delta^{i_2}}{\longrightarrow} [q+2] \longrightarrow \
\ldots \ \longrightarrow [n-1] \overset{\delta^{i_{n-q}}}{\longrightarrow} [n],
\end{split}
\end{equation*}
and the factorization of $\al$ is well-defined up to cosimplicial
identities.  Suppose that $\vec{a} \in A^{\ot (n+1)}$ and
$$ g: C( \Delta [n]) \to k[U_{*+1}]  $$
is a chain map.  We define
$$  \psi (g \ot \vec{a}) (\al ) = g(\al) \ot (d_{i_1} \circ d_{i_2}
\circ \cdots \circ d_{i_{n-q}} (\vec{a})).  $$
Now, $H_*(D_*(k[u_{*+1}]) \ot N^{\rm cy}(A)) \simeq HH_*(A)$, since
$D_*(k[U_{*+1}])$ is an acyclic simplicial $k$-module.  Also,
$$H_*\big( D_*\big( k[u_{*+1}] \ot N^{\rm cy}(A)\big) \big) \simeq
  H_*\big( k[U_{*+1}] \ot N^{\rm cy}(A) \big) .  $$

Define an inclusion of chain complexes
\begin{equation*}
\begin{split}
&  N^{\rm cy}(A) \to k[U_{*+1}] \ot N^{\rm cy}(A)  \\
&  A^{\ot (n+1)} \to k[U_{n+1}] \ot A^{\ot (n+1)}  \\
&  \vec{a} \mapsto \tau_{n+1} \ot \vec{a},
\end{split}
\end{equation*}
where $\tau_{n+1}$ is the cyclic shift in $U_{n+1}$.  We then have a
commutative diagram of chain complexes, where the diagonal and
vertical arrows are inclusions:
$$ \CD
N^{\rm cy}(A) @>>> k[U_{*+1}] \ot N^{\rm cy}(A)  \\
@VVV \hskip-1in \swarrow \hskip1in @VVV   \\
D_*(k[U_{*+1}]) \ot N^{\rm cy}(A) @>>\psi > D_*\big(
k[U_{*+1}] \ot N^{\rm cy}(A) \big). \endCD  $$
On homology
\begin{equation} \label{5.1}
\CD
HH_*(A) @>>> H_*(k[U_{*+1}] \ot N^{\rm cy}(A)) \\
@V\simeq VV \hskip-1.1in \swarrow \hskip1.1in @VV\simeq V \\
H_*(D_*(k[U_{*+1}]) \ot N^{\rm cy}(A)) @>>\psi_* >
H_* \big( D_*\big( k[U_{*+1}] \ot N^{\rm cy}(A) \big)\big).  \endCD  
\end{equation}
It follows that 
$$  HH_*(A) \to H_*(k[U_{*+1}] \ot N^{\rm cy}(A))  $$
is an isomorphism, and $HH_*(A) \simeq H_*(P_*(A))$.  In fact, every
arrow in diagram \eqref{5.1} is an isomorphism.

\end{document}